\documentclass[amstex,12pt,reqno]{amsart}
\usepackage{amsmath, amssymb, amsthm}
\usepackage{mathrsfs}

\usepackage{mathtools}
\mathtoolsset{showonlyrefs=true}

\newtheorem{Theorem}{Theorem}[section]
\newtheorem{Proposition}[Theorem]{Proposition}

\newtheorem{Lemma}[Theorem]{Lemma}
\newtheorem{Corollary}[Theorem]{Corollary}
\theoremstyle{definition}
\newtheorem{Definition}[Theorem]{Definition}
\theoremstyle{remark}
\newtheorem{Example}[Theorem]{Example}

\newtheorem*{remark*}{Remark}

 \providecommand\CC{{\mathbb C}}

\providecommand\brs{\begin{remark*}}
\providecommand\ers{\end{remark*}}

\providecommand\be{\begin{enumerate}}
\providecommand\ee{\end{enumerate}}
\providecommand\bT{\begin{Theorem}}
\providecommand\eT{\end{Theorem}}
\providecommand\bP{\begin{Proposition}}
\providecommand\eP{\end{Proposition}}
\providecommand\bD{\begin{Definition}}
\providecommand\eD{\end{Definition}}
\providecommand\bE{\begin{Example}}
\providecommand\eE{\end{Example}}
 \providecommand\bL{\begin{Lemma}}
\providecommand\eL{\end{Lemma}}
\providecommand\bC{\begin{Corollary}}
\providecommand\eC{\end{Corollary}}

\providecommand\bpp{\begin{proof}} \providecommand\epp{\end{proof}}
\providecommand\bee{\begin{equation}}
\providecommand\eee{\end{equation}}

  \providecommand\beqq{\begin{eqnarray*}}
\providecommand\eeqq{\end{eqnarray*}}

\providecommand\besp{\begin{split}}
\providecommand\eesp{\end{split}}
\providecommand\bay{\begin{array}}
\providecommand\eay{\end{array}}
\providecommand\sgn{\operatorname{sgn}}
\providecommand\bes{\begin{equation}\begin{split}}

\def\f#1#2{\frac{#1}{#2}}

\def\<{\leq}
\begin{document}

\title[Toeplitz operator ]{A counterexample of the Fredholm of  Toeplitz operator}
\author{Liu Hua}
\author{Zhang Xinyang}

\begin{abstract} In this paper we study the essential spectra of the Toeplitz operator on the Hardy space $H^1$. We give a counterexample to show that the Toeplitz operator with symbol is not Fridholm, which give a counterexample to the conjecture by J.A. Virtanen J A in 2006.

\end{abstract}

\keywords{ Hardy space,  Toeplitz operator, continuous symbol,  $\textbf{BMO}_{\log}$}
\subjclass[2000]{Primary: 42B20,43A65, 44A15; Secondary: 22D30,  30E25, 30G35, 32A35}

\address{hualiu@tute.edu.cn, Department of Mathematics, Tianjin
University of Technology and Education, Tianjin 300222, China}
\address{zhangxystd@163.com, School of Mathematics and Computer Science, Ningxia Normal University, Gu Yuan756000, China}

\maketitle

\bigskip

\section{Introduction}
Denote by $\mathbb{D}$ the unit disc in $\CC$  and $\mathbb{T}$ its boundary. Let $f(t)\in L^p(p\geq1)$. Define
$$
Pf(z)=\f1{2\pi i}\int_{\mathbb{T}}\f{f(t}{t-z}dt,\
\ |z|<1,
$$
and
$$
Qf(z)=-\f1{2\pi i}\int_{\mathbb{T}}\f{f(t}{t-z}dt,\  |z|>1.
$$
For $t\in \mathbb{T}$, denote by $Pf(t)=\lim_{z\rightarrow t,z\in \mathbb{D}}Pf(z)$ and $Qf(t)=\lim_{z\rightarrow t,z\in \mathbb{D}}Qf(z)$.
Denote by $H^p(\mathbb{T})$ the Hardy space and $H_0^p(\mathbb{T})$ the collections of the Hardy functions which vanishing at $0$. In the following, we prefer to $\overline{H_0^p(\mathbb{T})}$  The Toeplitz operator ${T}_a$ with the symbol $a\in L^\infty$ is defined by
\begin{equation}\label{top}
  {T}_a:\ H^p(\mathbb{T})\rightarrow H^p(\mathbb{T}),\ f\mapsto P(af).
\end{equation}

When $1<p<\infty$ and $a\in L^{\infty}(\mathbb{T})$,  the Toeplitz operator $T_a$ is bounded on $H^p (\mathbb{T})$. Moreover, if $a$ is continuous, the Fredholm properties of $T_a$ is well known: $T_a$ is Fredholm if and only $a$ is nowhere zero({\it see} \cite{bs}). When $p=1$, it is much more complex. The Toeplitz operator $T_a$ with
symbol $a$ is bounded on $H^1$ if and only if $a\in L^\infty$  and $Qa \in\mathbf{BMO}_{\log}$, \cite{sw,jp,to}. But it is more difficult to consider the Fredhlolmness of $T_a$ on $H^1$, \cite{pv,v}. In \cite{v}, J.A. Vertanen conjectured that the Fredholmness of $T_{a}$ still holds for continuous symbol  $a$ satisfying $Qa\in BMO_{\log}$.

In this paper, we give an example $a(z)=(1-\ln(1+z))/(1-\ln(1+z^{-1}))$ to disapprove the Vertanen's conjecture. We give the outline of the proof as follows.

 Notice that $a$ is  continuous on $\mathbb{T}$ since $a-1$ and $\ln a$ are equivalent infinitesimal when $z$ tends to $-1$ on $\mathbb{T}$.
  Obviously, $\textbf{ind}a=0$, here $\textbf{ind}a$ is the index of $a$ is continuous on $\mathbb{T}$. First we will prove that $\ln a^{+}$, $\ln a^{-}$ satisfy the assumption of the Vertanen's conjecture, i.e., both $Pa$
and $Qa$ are $\textbf{BMO}_{\log}$ functions, by approximate their derivative.
Then we show that the kernel of $T_{a}$ is $0$ and its image is dense in (i.e. closure is) $H^{1}(\mathbb{T})$. Finally, we construct
 a function in $H^{1}(\mathbb{T})$ which has no preimage of $\mathbb{T}_\alpha$ in $H^{1}(\mathbb{T})$. These two facts indicate  that the image of $T_a$ is closed in  $H^{1}(\mathbb{T})$, i.e, $T_a$ is not Fredholm.
 
 \noindent The second author is the main contributor to this paper.
\section{Preliminary results}

In this article we choose the continuous branch of the argument function $\arg$ such that $\arg(z) \in (-\pi, \pi]$ for $z\in \mathbb{C}$. And the the logarithm function takes is defined by {$\textbf{Im} \ln(z)=\arg(z)$ on $\mathbf{C}\backslash\{z|\textbf{Re}z\leq 0, \textbf{Im}z=0\}$. We denote $\theta=\theta(z)=\arg(z)$ for $z\neq -\left\vert z\right\vert$ and $\vartheta=\vartheta(z)=\arg(-z)$ for $z\neq \left\vert z\right\vert$.

\bE We consider a function $a=a^{+}/a^{-}$ on the unit circle $\mathbb{T}$, in which $a^{+}(z)=1-\ln(1+z)$ analytic in
$\mathbf{C}\backslash\{z|\textbf{Re}z\leq -1, \textbf{Im}z=0\}$ and $a^{-}(z)=1-\ln(1+z^{-1})$ analytic in $\mathbf{C}\backslash\{z|-1\leq \textbf{Re}z\leq 0, \textbf{Im} z=0\}$. Since $\textbf{Re}a^{+}(z)=1-\ln|1+z|>1-\ln2>0$ for $z\in \mathbf{D}\cup \mathbb{T}\backslash\{-1\}$, while
$\textbf{Re}a^{-}(t)=1-\ln|1+t^{-1}|>1-\ln2>0$ for $t\in \mathbb{T}\backslash\{-1\}$. Then
 we have that $\ln (t)=\ln a^{+}(t)-\ln a^{-}(t)$ holds on $\mathbb{T}\backslash\{-1\}$ in the sense of main value, which is a pure imagine odd function of $\theta \in (-\pi, \pi)$. When $z=-\exp(i\vartheta)$ tends to $-1$ on $\mathbb{T}$, $\textbf{Re}\ln(1+z)=\ln(2\sin(\left\vert \vartheta\right\vert/2))$ which is equivalent infinity of $\ln\left\vert \vartheta\right\vert$, while $\textbf{Im}\ln(1+z)=\theta/2$ tends to $\pm\pi/2$ respectively when $t$ tends to $-1$ from above or below.   Thus both $a^{+}$ and $  a^{-} $ are belong to $ L^{1}(\mathbb{T})$.

 Since $\ln(1+t)$ has only weak singularity at $-1$, we have that
 $$\ln(1+z)=\int_{\mathbb{T}}\f{\ln(1+t)}{t-z}dt,\ z\in \mathbb{D},$$
 i.e. $a^{+}\in H^{1}(\mathbb{T})$. Similar argument gives  $a^{-}\in H^1{\mathbb{C}\setminus \mathbb{D}}\subset \overline{H^{1}(\mathbb{T})}$. Moreover when $t\in \mathbf{D}$ tends to $-1$ from above or below, $\textbf{Re}\ln a^{+}=\textbf{Re}\ln a^{-}$ are equivalent to $\ln(-\ln\left\vert \vartheta\right\vert )$, while $\textbf{Im}\ln a^{+}=-\textbf{Im}\ln a^{-}$ are equivalent to $\pm\pi/2\ln(\left\vert \vartheta\right\vert)$.  $\ln a^{+}(t)$ still has anly weak singularity at $-1$. We again obtain that $\ln a^{+}$ belongs to $H^{1}(\mathbb{T})$ and $\ln a^{-}$  belongs to $\overline{H^{1}(\mathbb{T})}$, respectively.

 Furthermore, $1/a^{+}\in H^{\infty}(\mathbb{T})$ and $1/a^{-}\in \overline{H^{\infty}(\mathbb{T})}$.
 All of them are  continuous at $\mathbb{T}\backslash\{-1\}$ with the logarithmic singularity at $-1$.
  Moreover $\ln a(t)$ tends to $0$  and equivalent to $\pm\frac{\sgn(\pi}{\ln(\pi -\left\vert \theta \right\vert)}$ when $t$ tends to $-1$ from above or below.
\eE

Let the connecting arc $I\subset \mathbb{T}$ (we also call  $I$ by the interval). Without confusion, we assume $-1\notin I$ and identity $I$ with the real set $\{\theta\in [0,2\pi),e^{i\theta}\in I\}$.
Let $f(z)\in H^1$, denote by
$$
f_I=\f1{|I|}\int_I|f(e^{i\theta})|d\theta
$$
and
$$
\mathrm{MO}_I(f)=\f1{|I|}\int_I|f(e^{i\theta})-f_I|d\theta.
$$
\bE Let $I=[0,d]$.
\begin{eqnarray}
 \nonumber \mathrm{MO}_I(\ln x)&=& \f1{d}\int_I |\ln \vartheta-\f1d\int_I \ln \xi d\xi|dx\\
   &=&\f1d \int_I|\ln x- \ln d+1|dx=2e^{-1}.
\end{eqnarray}
\eE
\noindent In this paper we need the following lemma.
\bL
Let $f(x)$ be a bounded variation function on $I=[0,d]$ with the total variation $V_I(f)$. We obtain.
\begin{equation}\label{b-v}
  \mathrm{Mo}_I(f)\<\f12V_I(f).
\end{equation}
\eL
\bpp
Let $I=[a,b]\subset [0,d]$ and $f(x)\in C^1[a,b]$. We obtain
\begin{eqnarray}\label{le}
 \nonumber \int^b_a\left|f(x)-f_I\right|dx  &=& \int^b_a\left|f(x)-\f1{|I|}\int^b_af(t)dt\right|dx \\
   \nonumber &=&  \f1{|I|}\int^b_a\left|(\int^x_a+\int^b_x)(\int^x_t f'(s)ds)dt\right|dx\\
 \nonumber   &\leq& \f1{|I|}\int^b_a\left|\int^x_a\int^x_t f'(s)dsdt\right|dx+ \f1{|I|}\int^b_a\left|\int^b_x\int^x_t f'(s)dsdt\right|dx\\
 &=& \mathrm{I}+\mathrm{II}.
\end{eqnarray}
Notice that
\begin{eqnarray}
 \nonumber  \left|\int^x_a\int^x_t f'(s)dsdt\right|&\leq& \int^x_a\int^x_t |f'(s)|dsdt \\
 \nonumber   &=& \int^x_a(s-a)|f'(s)|ds,
\end{eqnarray}
 we have
\begin{eqnarray}\label{l1}
 \nonumber  \mathrm{I}&\leq&\f1{|I|}\int^b_a \int^x_a(s-a)|f'(s)|dsdx\\
   \nonumber  &=&\f1{|I|}\int^b_a( \int^b_s(s-a)|f'(s)|dx)ds\\
   &=&\f1{|I|}\int^b_a(s-a)(b-s)|f'(s)|ds.
\end{eqnarray}
Similarly, we also have
\begin{equation}\label{l2}
  \mathrm{II}\leq\f1{|I|}\int^b_a(s-a)(b-s)|f'(s)|ds.
\end{equation}

By (\ref{ll}), (\ref{l1} and (\ref{l2}, we obtain
\begin{eqnarray}\label{l3}
  \nonumber \int^b_a\left|f(x)-f_I\right|dx  &\leq& \f1{|I|}\int^b_a2(s-a)(b-s)|f'(s)|ds \\
   \nonumber &\leq& \f1{|I|} \int^b_a2(s-a)(b-s)|f'(s)|ds\\
    &\leq&\f1{|I|}\int^b_a\f{(s-a+b-s)^2}2|f'(s)|ds\\
    &=& \f{|I|}2\int^b_a|f'(s)|ds\leq\f{|I|}2V_I(f),
\end{eqnarray}
which gives (\ref{b-v}).

\epp
\noindent \bf{Remark.} \rm We only proved (\ref{b-v}) for $f\in C^1[a,b]$ because it is enough for this article. In general, let $f$ be a bounded variation function. Since $f(x)$ is Riemann integrable, apply the Abel summation formula to the Riemann summation of the concerned integral, we can complete the proof of the lemma by Riemann summation and Ablel formula.

\bD $^\text{\cite{pv}}$ A function $f$ is in $\mathrm{BMO}_{\log}$ if $f\in L^1(\mathbb{T})$ and
\begin{equation}\label{bmolog}
  ||f||_{\log}=\sup_I\f{4\pi}{\log{|I|}}\mathrm{MO}_I(f)<+\infty,
\end{equation}
where the supremum is taken over all arcs $I$ of $\mathbb{T}$. The space $\mathrm{BMO}_{\log}$ is a Banach space under the norm $||f||_{\mathrm{BMO}_{\log}}=\f1{2\pi}|\int f(e^{i\theta})d\theta|+ ||f||_{\log}$.
\eD
In the next section, We will show that $\ln(-\ln(\left\vert x \right\vert))$ and $\sgn(x)/\ln(\left\vert x \right\vert)$ are $\textbf{BMO}_{\log}$ in the neighbourhood of $x=0$.  $\ln a^{+}$ and $\ln a^{-}$ seem to be $\textbf{BMO}_{\log}$ as they are equivalent infinity of $\ln(-\ln\left\vert x \right\vert)$. However, equivalent infinity is not enough to keep the $\textbf{BMO}_{\log}$ property. In order to prove it we should estimate their variety and conjecture $\textbf{BMO}_{\log}$ property keeps if the variety of a function is controlled by a $\textrm{BMO}_{\log}$ function everywhere.
\section{$\mathrm{BMO}_{\log} space$}
\bL \label{lemma1}
If real functions $f$ and $g$ on interval $I=[a,b]$ satisfy $0\leq \left\vert f(x_{1})-f(x_{2})\right\vert \leq \left\vert g(x_{1})-g(x_{2})\right\vert$ for all pair of points $x_{1}, x_{2}\in I$, then $MO_{I}(f)\leq 2MO_{I}(g)$. Moreover, if $f$ and $g$ are monotaneous, then $MO_{I}(f)\leq MO_{I}(g)$.\eL
\bpp
Consider $y^{+}=\inf\{g(x)|x\in I\ {\rm and}\  g(x)\geq g_{I}\}$, and $y^{-}=\sup\{g(x)|x\in I\  {\rm and}\  g(x)\leq g_{I}\}$, let $0\leq \lambda \leq 1$ such that $\lambda y^{+}+(1-\lambda)y^{-}=g_{I}$, and $x_{1},x_{2}\in I$ such that both $0\leq g(x_{1})-y^{+}\leq \epsilon$ and $0\leq y^{-}-g(x_{2})\leq \epsilon$ hold. Denote by $y=\lambda f(x_{1})+(1-\lambda)f(x_{2})$. Notice that $\int_{I}( f(x)-f_{I}) dx=0$, if $y\leq f_{I}$ we have
\begin{eqnarray}
 \nonumber\int_{I}\left\vert f(x)-f_{I} \right\vert dx  &=& 2\int_{I}(f(x)-f_{I})^{+}dx \\
 \nonumber  &\leq& 2\int_{I}(f(x)-y)^{+}d \\
   &\leq& 2\int_{I}(\lambda(f(x)-f(x_{1}))^{+}+(1-\lambda)(f(x)-f(x_{2}))^{+})dx
\end{eqnarray}

or if $y\geq f_{I}$ we have
\begin{eqnarray}
 \nonumber\int_{I}\left\vert f(x)-f_{I} \right\vert dx  &=& 2\int_{I}(f(x)-f_{I})^{-}dx \\
 \nonumber  &\leq& 2\int_{I}(f(x)-y)^{-}d \\
   &\leq& 2\int_{I}(\lambda(f(x)-f(x_{1}))^{-}+(1-\lambda)(f(x)-f(x_{2}))^{-})dx.
\end{eqnarray}

In the above either cases we obtain
\begin{eqnarray}
 \nonumber\int_{I}\left\vert f(x)-f_{I} \right\vert dx &
\leq & 2\int_{I}(\lambda\left\vert f(x)-f(x_{1})\right\vert +(1-\lambda)\left\vert f(x)-f(x_{2})\right\vert )dx\\
\nonumber&\leq & 2\int_{I}(\lambda\left\vert g(x)-g(x_{1}) \right\vert +(1-\lambda)\left\vert g(x)-g(x_{2})\right\vert )dx\\
\nonumber&\leq &2\int_{I}(\lambda\left\vert g(x)-y^{+}\right\vert +(1-\lambda)\left\vert g(x)-y^{-}\right\vert +\epsilon)dx\\
&=&2\int_{I}(\left\vert g(x)-\lambda y^{+}-(1-\lambda)y^{-}\right\vert +\epsilon)dx.
\end{eqnarray}
The last identity is due to that either $g(x)\geq y^{+}$ or $g(x)\leq y^{-}$. In short, $MO_{I}(f)\leq 2(MO_{I}(g)+\epsilon)$.

When $f$ and $g$ are monotaneously increasing functions.  Then there exists $x_0$ such that
\begin{eqnarray}
 \nonumber&\ & \int_{I}(\lambda(f(x)-f(x_{1}))^{+}+(1-\lambda)(f(x)-f(x_{2}))^{+})dx\\
\nonumber&\leq & \int_{I}(\lambda(g(x)-g(x_{1}))^{+}+(1-\lambda)(g(x)-g(x_{2}))^{+})dx\\
&\leq & \int_{I}(\lambda(g(x)-y^{+})^{+}+(1-\lambda)(g(x)-y^{-})^{+}+\epsilon)dx,
\end{eqnarray}
and
\begin{eqnarray}
 \nonumber&\ & \int_{I}(\lambda(f(x)-f(x_{1}))^{-}+(1-\lambda)(f(x)-f(x_{2}))^{-})dx\\
\nonumber&\leq & \int_{I}(\lambda(g(x)-g(x_{1}))^{-}+(1-\lambda)(g(x)-g(x_{2}))^{-})dx\\
&\leq & \int_{I}(\lambda(g(x)-y^{+})^{-}+(1-\lambda)(g(x)-y^{-})^{-}+\epsilon)dx.
\end{eqnarray}
Similarly, we obtain
\begin{eqnarray}
 \nonumber&\ & \int_{I}(\lambda(g(x)-y^{+})^{+}+(1-\lambda)(g(x)-y^{-})^{+})dx\\
\nonumber&= & \int_{I}(\lambda(g(x)-y^{+})^{-}+(1-\lambda)(g(x)-y^{-})^{-})dx\\
&= & \int_{I}(\left\vert g(x)-\lambda y^{+}-(1-\lambda)y^{-}\right\vert )dx/2.
\end{eqnarray}
Another proof is $x_{0}\in I$ such that $g(x_{0})\leq g_{I}$ for $I_{2}:x<x_{0}$ and $g(x_{0})\geq g_{I}$ for $I_{1}:x<x_{0}$, hence $y^{+}$ and $y^{-}$ are just the right and left limit of $g$ at $x_{0}$, respectively. Take $x_{1}=x_{0}^{+}$ and $x_{2}=x_{0}^{-}$. Denote by $I_1=\{x\in I,x<x_0\}$ and $I_2=\{x\in I,x<x_0\}$. It is easy to check that $f(x)\leq g(x)-g(a^+)+f(a^+)$ and hence $f_I\<g_I$. Then we have
\begin{eqnarray}
   \nonumber\int_{I}(f(x)-y)^{+}dx&=& \int_{I_{1}}(y-f(x))dx \\
  \nonumber &=& \int_{I_{1}}(\lambda (f(x_0^+)-f(x))+(1-\lambda) (f(x_0^-)-f(x))dx\\
  \nonumber &\leq& \int_{I_{1}}(\lambda (f(x_0^+)-f(x))+(1-\lambda) (f(x_0^-)-f(x))dx\\
  \nonumber &\leq&\int_{I_{1}}(\lambda (g(x_0^+)-g(x))+(1-\lambda) (g(x_0^-)-f(x))dx\\
   &=&\int_{I_{1}}(g_I-g(x))dx
\end{eqnarray}
$$
\int_{I}(f(x)-y)^{+}dx=\int_{I_{1}}(f(x)-y)dx
$$ and $\int_{I}(f(x)-y)^{-}dx=\int_{I_{2}}(f(x)-y)dx$.
\epp

\bC\label{cor1}
If real  functions $f$ and $g$ on interval $I = [a,b]$ satisfy $\|f'(x)\|\leq g'(x)$, then $MO_{I}(f)\leq 2MO_{I}(g)$.
\eC
\bpp
It is easy to check by Lemma \ref{lemma1} since
$$
|f(x_1)-f(x_2)|=|\int_{x_1}^{x_2} f'(x)dx|\leq\int_{x_1}^{x_2} |f'(x)|dx\leq\int_{x_1}^{x_2} g'(x)dx=g(x_2)-g(x_1)$$
for $x_1,x_2\in I$ and $x_1<x_2$.\epp

\bP
If even real function $f(x)$ satisfies $\left\vert f'(x) \right\vert \leq -\left\vert x \right\vert ^{-1}\ln^{-1}\left\vert x \right\vert$
when $0<\left\vert x \right\vert \leq \delta\leq\f14 \exp(-1)$,
then $MO_{I}(f)\leq -8e^{-1}\ln^{-1}(\left\vert I \right\vert/2)$ for every interval $I\subseteq [-\delta,\delta]$. Moreover if $f(x)$ is monotaneous on $0<x\leq \delta$ then the 8 above can be reduced to 4.\eP
\bpp
Let $I=[\delta_{l},\delta_{r}]$. Recall that $-x^{-1}\ln^{-1} x$ is decreasing for $0<x<e^{-1}$. If $\delta_{l}>0$ then

\begin{eqnarray}\label{d-com}
 \nonumber\left\vert f'(x+\delta_{l})\right\vert &\leq & -(x+\delta_{l})^{-1}\ln^{-1}(x+\delta_{l})\\
\nonumber&\leq & -x^{-1}\ln^{-1}(x)\\
&\leq & -\ln^{-1}\left\vert I\right\vert x^{-1}=(-\ln^{-1}\left\vert I\right\vert \ln x)'
\end{eqnarray}
for $x\in [0,\left\vert I \right\vert ]$ since $x+\delta_l<e^{-1}$. By and (\ref{d-com}) Corollary \ref{cor1}, we obtain
\begin{eqnarray}
  \nonumber\mathrm{MO}_{I}(f)&= & \mathrm{MO}_{[0,\left\vert I \right\vert ]}(f(x+\delta_{l}))\\
 \nonumber&\leq & -2\ln^{-1}(\left\vert I \right\vert )\mathrm{MO}_{[0,\left\vert I \right\vert ]}(\ln (x))=-4e^{-1}\ln^{-1}(\left\vert I\right\vert )\\
&=&-8e^{-1}\f1{2\ln |I|}\leq  -8e^{-1}\ln^{-1}(\left\vert I\right\vert/2).
\end{eqnarray}

In the case $-\delta_{r}\leq \delta_{l}<0$, we define the  piecewise function $g(x)$ on $[0,\left\vert I \right\vert ]$ valued $f(x/2)$ when $x\leq -2\delta_{l}$ and $f(x+\delta_{l})$ when $x\geq -2\delta_{l}$. Now
\[
g'(x)\leq -x^{-1}\ln^{-1}(x/2)\leq -\ln^{-1}(\left\vert I \right\vert/2)x^{-1}
\]
when $x<-2\delta_l$, and
\[
g'(x)\leq -(x+\delta_{l})^{-1}\ln^{-1}(x+\delta_{l})\leq -2x^{-1}\ln^{-1}(x/2)\leq -2\ln^{-1}(\left\vert I \right\vert/2)x^{-1}
\]
when $x>-2\delta_l$.
Notice that
\begin{eqnarray}
   \nonumber\int_{\delta_l}^{\delta_r} f(x)dx&=&\int_{\delta_l}^{-\delta_l}f(x)dx+\int_{-\delta_l}^{\delta_r}f(x)dx\\
   \nonumber&=&2\int_0^{-\delta_l}f(x)dx+\int_{-2\delta_l}^{|I|}f(x+\delta_l)dx\\
   &=&\int_0^{-2\delta_l}g(x)dx+\int_{-2\delta_l}^{|I|}g(x)dx=\int_0^{|I|}g(x)dx.
\end{eqnarray}
We obtin
\begin{eqnarray}
\mathrm{MO}_{I}(f)
\nonumber&= & \mathrm{MO}_{[0,\left\vert I \right\vert ]}(g)\\
\nonumber&\leq & -4\ln^{-1}(\left\vert I \right\vert/2)\mathrm{MO}_{[0,\left\vert I \right\vert ]}(\ln (x))\\
&= & -8e^{-1}\ln^{-1}(\left\vert I \right\vert/2).
\end{eqnarray}
\epp
\bP
If function (not required real) $f(x)$ is continuous at 0 and satisfies $\left\vert f'(x) \right\vert \leq \left\vert x\right\vert ^{-1}\ln^{-2}\left\vert x\right\vert$ when $0<\left\vert x\right\vert \leq \delta \leq \exp(-2)$, then $MO_{I}(f)\leq \ln^{-1}(-\ln(\left\vert I\right\vert/2))$ for all $I\subseteq [-\delta,\delta]$. \eP
\bpp
Denote by $V_{I}(f)$ the total variation of $f(x)$ on $I$, then $MO_{I}(f)\leq V_{I}(f)/2$. Since $f$ is absolutely continuous on $I$, then $\mathrm{V}_{I}(f)=\int_{I}\left\vert f'(x)\right\vert dx$ and we have
\begin{eqnarray}
\nonumber \mathrm{MO}_{I}(f)&\<& \mathrm{V}_I(f)=\int_{I}\left\vert f'(x)\right\vert dx\\
\nonumber&\leq & \int_{I}\left\vert x\right\vert ^{-1}\ln^{-2}(\left\vert x\right\vert )dx/2\\
\nonumber&\leq & \int_{[0,\left\vert I\right\vert/2]}\left\vert x\right\vert ^{-1}\ln^{-2}(\left\vert x\right\vert )dx\\
&= & \ln^{-1}(-\ln(\left\vert I\right\vert/2)),
\end{eqnarray}
the last unequal holds because comparing $[0,\left\vert I\right\vert/2]$ with $[\delta_{l},\delta_{l}+\left\vert I\right\vert/2]$ and $[\delta_{l}+\left\vert I\right\vert/2,\delta_{l}+\left\vert I\right\vert ]$ when $\delta_{l}>0$, or  $[-\delta_{l},\left\vert I\right\vert/2]$ with $[\left\vert I\right\vert/2,\delta_{r}]$ when $-\delta_{r}\leq \delta_{l}<0$.
\epp
\bT $a $  is $BMO_{\log}$ on $\mathbb{T}$\eT
\bpp
According to the two propositions above, the problem of checking $BMO_{\log}$ property is converted to comparing the derivatives with $-\left\vert x\right\vert ^{-1}\ln^{-1}(\left\vert x\right\vert )$ or $\left\vert x\right\vert ^{-1}\ln^{-2}(\left\vert x\right\vert )$.
\vskip2mm
\noindent It is easy to check that
\begin{eqnarray}\label{d-lna}
\nonumber \f{ d\ln(a^{+}(-e^{ix}))}{dx}&=&-i(1-e^{-ix})^{-1}(1-\ln(1-e^{(ix}))^{-1}\\
& \sim&
 x^{-1}\ln^{-1}\left\vert x\right\vert
\end{eqnarray}
 as $x\rightarrow0$. Moreover we have
\begin{eqnarray}\label{d-ilna}
  \nonumber \f{d\textbf{Im}\ln(a^{+}(-e^{ix}))}{dx} &=& (4\sin^{-2}(\left\vert x\right\vert/2)((1-\ln(2\sin(\left\vert x\right\vert/2)))^2 \\
  \nonumber &\ & +(\pi/2-\left\vert x\right\vert/2)^2)^{-1}((\cos(x)-1)(1-\ln(2\sin(\left\vert x\right\vert/2))) \\
 \nonumber  &\ &  +(-\sin(x))(x-\sgn(x)\pi)/2)\\
 &\sim& (\pi/2)\left\vert x\right\vert ^{-1}\ln^{-2}\left\vert x\right\vert
\end{eqnarray}
 as $x\rightarrow0$. By (\ref{d-lna}) and (\ref{d-ilna}) we obtain
 \begin{equation}\label{d-lnaa}
\f{ d\textbf{Re}\ln(a^{+}(-e^{ix}))}{dx}\sim
 x^{-1}\ln^{-1}\left\vert x\right\vert
 \end{equation}
  as $x\rightarrow0$.

Take an arbitray $\delta$ small enough that $\left\vert dRe\ln(a^{+}(-\exp(ix)))/dx\right\vert \leq -M_{0}\left\vert x\right\vert ^{-1}\ln^{-1}\left\vert x\right\vert$, denote $M_{1}=\max_{\left\vert \theta \right\vert \leq \pi-\delta/2}\left\vert dRe\ln(a^{+}(\exp(i\theta)))/d\theta \right\vert$, for arc $I\subseteq T$ if $\left\vert I\right\vert \geq \delta/2$ then $MO_{I}(Re\ln(a^{+}))\leq 2\int_{-\pi}^{\pi}\left\vert \ln(a^{+}(\exp(i\theta)))\right\vert d\theta/(\delta/2)$, if $I\subseteq \{\exp(i\theta):\delta/2-\pi \leq \theta \leq \pi-\delta/2\}$ then $MO_{I}(Re\ln(a^{+}))\leq M_{1}\left\vert I\right\vert/2$, if $I\subseteq \{-\exp(ix)|-\delta \leq x \leq \delta\}$ then $MO_{I}(Re\ln(a^{+}))\leq -4\exp(-1)M_{0}\ln{-1}(\left\vert I\right\vert/2)$, thus $Re\ln(a^{+})=Re\ln(a^{-})$ is $BMO_{\log}$ on $\mathbb{T}$ and $Im\ln(a^{+})=-Im\ln(a^{-})=\ln(a)/2i$ is similar, and so does $a$ since
\begin{eqnarray}
 \nonumber
 V_{I}(a)
&= & \int_{I}\left\vert da(-e^{ix})/dx\right\vert dx\\
 \nonumber &= & \int_{I}\left\vert a(-e^{ix})d\ln(a(-e^{ix}))/dx\right\vert dx\\
 \nonumber& = & \int_{I}\left\vert d\ln(a(-e^{ix}))/dx\right\vert dx\\
&= & V_{I}(\ln (a)).
\end{eqnarray}
\epp
\vskip2mm
\noindent Now we will prove that $Qa$ is $BMO_{\log}$, for which we first give a simple representation of $Qa$. Denote by $\tilde{f}(z)=f(z^{-1})$. We have
 \begin{equation}\label{}
   (Q\tilde{f})(z)=-\f1{2\pi i}\int_{\mathbb{T}}\f{f(z^{-1})}{t-z}dt=\widetilde{Pf}(z)+f_0\f1{2\pi i}\int_{\mathbb{T}}\f{f(t)}tdt+f_0,
 \end{equation}
 where $f_0=\f1{2\pi i}\int_{\mathbb{T}}\f{f(t)}tdt$. Thus we obtain
 \begin{equation}\label{gggg}
   Qa(z)=a-Pa(z)=a-\widetilde{\widetilde{Pa}}(z)=a-\widetilde{Q\tilde{a}}(z)+a_0.
 \end{equation}
 Recall that
 $\tilde{a^{+}}=a^{-}=1-\ln(1+z^{-1})$ is also analytic in $\mathbf{D}\backslash \{z|-1<Rez\leq0,Imz=0\}$, splitting into two branches at $\{z|-1<Rez<0,Imz=0\}$ and has weak (logarithmic order) singularity near $0$. Denote by $L^{\pm}_{1-\epsilon}$ respectively the upper  and lower sides of the interval $[-1,-\epsilon]$, $L_\epsilon$ the circle with the radius of $\epsilon$ centered on $O$, which directions are set to correspond to the part of the boundary of $\mathbf{D}\backslash \{z|-1<Rez\leq0,Imz=0\}$.

 By Cauchy theorem we obtain
\[
\int_{\mathbb{T}+L^+_{1-\epsilon}+L_{\epsilon}+L^-_{l-\epsilon}}\frac{1-\ln(1+t^{-1})}{(t-z)(1-\ln(1+t))}dt=0.
\]
It is easy to check that
\begin{equation}\label{l-ep1}
\int_{L^+_{1-\epsilon}+L^-_{1-\epsilon}}\frac{1-\ln(1+t^{-1})}{(t-z)(1-\ln(1+t))}dt=\int_{-1}^{\epsilon}\frac{2\pi i}{(t-z)(1-\ln(1+t))}dt
\end{equation}
and
\begin{equation}\label{l-ep2}
\left\vert\int_{L_{\epsilon}}\frac{1-\ln(1+t^{-1})}{(t-z)(1-\ln(1+t))}dt\right\vert \leq 2\pi\epsilon\frac{\ln(1+\epsilon^{-1})+\pi -1}{(\left\vert z\right\vert -\epsilon)(1-\ln(1-\epsilon))}.
\end{equation}
Recall that (\ref{l-ep2}) tends to $0$ as  $\epsilon\rightarrow0$. We have
\begin{eqnarray}\label{q}
\nonumber Q\tilde{a}(z)&=& -\f1{2\pi i}\int_{\mathbb{T}}\frac{1-\ln(1+t^{-1})}{(t-z)(1-\ln(1+t))}dt\\
\nonumber&=&\f1{2\pi i}\int_{L^+_{1-\epsilon}+L_\epsilon+L^-_{1-\epsilon}}\frac{1-\ln(1+t^{-1})}{(t-z)(1-\ln(1+t))}dt \\
\nonumber&=&\int_{-1}^{\epsilon}\frac{2\pi i}{(t-z)(1-\ln(1+t))}dt+o(1)\\
 &=&\int_{-1}^{0}\frac{dt}{(t-z)(1-\ln(1+t))}.
\end{eqnarray}
Hence we obtain
\begin{eqnarray}\label{dq}
\nonumber\f{d}{dz}Q\tilde{a}(z)&= & \int_{-1}^{0}\frac{dt}{(t-z)^{2}(1-\ln(1+t))}\\
&= &  z^{-1}+\int_{-1}^{0}\frac{dt}{(t-z)(1+t)(1-\ln(1+t))^{2}}
\end{eqnarray}

Let $\theta=\arg(-z)$ for $|z|=1$. So $i\theta\sim-1-z$ as $z$ tends to $-1$. By (\ref{q}), it is easy to check that  $ReQ\tilde{a}(z)$ and $ImQ\tilde{a}(z)$ are even and odd respectively of $\theta$. So we assume that $\theta>0$.

To estimate the order of the integral in (\ref{dq}) as $z$ tends to $-1$ on $\mathbb{T}$,
 we divide it into three parts as follows
\begin{eqnarray}\label{3parts}
\nonumber\int_{-1}^{0}q(t,e^{i\theta})dt&=&\left\{\int_{-1}^{\theta-1}+\int_{\theta-1}^{e^{-3}-1}+\int_{e^{-3}-1}^{0}\right\}q(t,e^{i\theta})dt\\
&=&I_1(z)+I_2(z)+I_3(z),
\end{eqnarray}
where $q(t,z)$ is the integrand in (\ref{dq}).

Rewrite $I_1(z)$ as follows
\begin{eqnarray}\label{i1}
 \nonumber I_1(z)&= &(1+z)^{-1}\int_{-1}^{0}(\frac{1}{t-z}-\frac{1}{1+t})\frac{dt}{(1-\ln(1+t))^{2}}\\
 &=&(1+z)^{-1}I_{11}-(1+z)^{-1}I_{12}.
\end{eqnarray}

It is easy to check that
\begin{equation}\label{i12}
I_{12}=\int_{-1}^{\theta-1}(1+t)^{-1}(1-\ln(1+t))^{-2}dt=(1-\ln\theta)^{-1}.
\end{equation}

Since $\left\vert t-z\right\vert=|t-\cos\theta-i\sin\theta| \geq \sin\theta$ and $(1-\ln(1+t))^{-2}\leq (1-\ln\theta)^{-2}$, we have
\begin{eqnarray}\label{i11}
\nonumber|I_{11}|=\left\vert\int_{-1}^{\theta-1}(t-z)^{-1}(1-\ln(1+t))^{-2}dt\right\vert &\leq& \theta(\sin\theta)^{-1}(1-\ln\theta)^{-2}\\
&\sim&(1-\ln\theta)^{-2}.
\end{eqnarray}

By (\ref{i1})---(\ref{i11}) we obtain
\begin{equation}\label{I1}
  I_1(z)\sim(1+z)^{-1}I_{12}\sim(1+z)^{-1}\ln^{-1}\theta\sim(1+z)^{-1}\ln^{-1}(1+z).
\end{equation}

It is obvious that
\[
I_3(z)=\int_{e^{-3}-1}^{0}\frac{(1-\ln(1+t))^{-2}}{(t-z)(1+t)}dt
\]
is convergent as $z$ tends to $-1$, i.e.,
\begin{equation}\label{I3}
\lim_{z\rightarrow -1}I_3(z)=\int_{e^{-3}-1}^{0}(1+t)^{-2}(1-\ln(1+t))^{-2}dt.
\end{equation}

We still rewrite $I_2(z)$ as
\begin{equation}\label{i20}
 I_2(z)= (1+z)^{-1}\int_{\theta-1}^{e^-{3}-1}(\frac{1}{t-z}-\frac{1}{1+t})\frac{dt}{(1-\ln(1+t))^{2}}.
\end{equation}

Denote by $w=e^{i\varphi}$. To estimate $I_2(z)$ we consider
\begin{equation}\label{i21}
(1-\ln(t-w))^{-2}|^z_{-1}=-2i\int_{0}^{\theta}w(t-w)^{-1}(1-\ln(t-w))^{-3}d\varphi.
\end{equation}
Obviously there are $\left\vert t-w\right\vert \geq 1+t$ and $\left\vert 1-\ln(t-w)\right\vert \geq 1-\ln \left\vert t-w\right\vert$ for $\varphi$ small enough.

Let $|\varphi|<\f19$. We have
 $1+2(e^{-3}-1)\cos\varphi+(e^{-3}-1)^{2}\leq e^{-4}$. And then we obtain that $\left\vert t-e^{i\varphi}\right\vert \leq e^{-2}$ when $-1\leq t\leq e^{-3}-1$. Recall that  $s^{-1}(1-\ln s)^{-3}$ is decreasing when $s\in(0,\exp(-2))$. We obtain
\begin{equation}\label{i22}
\max_{0\leq x'\leq x}\left\vert t-z'\right\vert ^{-1}(1-\ln \left\vert t-z'\right\vert )^{-3}=(1+t)^{-1}(1-\ln(1+t))^{-3}.
\end{equation}
By (\ref{i22}) and
 \begin{equation}\label{i23}
 (1+t)^{-2}(1-\ln(1+t))^{-3}\leq -4\f d{dt}(1+t)^{-1}(1-\ln(1+t))^{-3}
 \end{equation}
 on $(-1,e^{-3}-1)$,
   we obtain
\begin{eqnarray}\label{i24}
\nonumber&\ & \left\vert \int_{\theta-1}^{e^{-3}-1}\f 1{t-z}((1-\ln(t-z))^{-2}-(1-\ln(t-w))^{-2})dt\right\vert \\
\nonumber&\leq & 2\theta\int_{\theta-1}^{e^{-3}-1}(1+t)^{-2}(1-\ln(1+t))^{-3}dt\\
\nonumber&\leq & -8\theta(1+t)^{-1}(1-\ln(1+t))^{-3}|^{e^{-3}-1}_{\theta-1}\\
&= & 8(1-\ln\theta)^{-3}-\f18e^3\theta.
\end{eqnarray}
Again recall that  $z^{-1}(1-\ln z)^{-2}$ is analytic on  $\{z|\theta-1+\cos\varphi\leq \textrm{Re}z\leq e^{-3}-1+\cos\varphi, \textrm{Im}z=\sin\varphi, 0\leq \varphi\leq \theta\}$. By Cauchy theorem we obtain
\begin{eqnarray}\label{i25}
\nonumber&\ &\int_{\theta-1}^{e^-{3}-1}(\frac{1}{t-z}-\frac{1}{1+t})\frac{dt}{(1-\ln(1+t))^{2}}\\
\nonumber&=& \int_{\theta-1}^{e^{-3}-1}(\f1{(t-z)(1-\ln(t-z))^2}-\f1{(t+1)(1-\ln(t+1))^2} )dt\\
&= & (\int_{e^{-3}}^{e^{-3}-z-1}-\int_{\theta}^{\theta-z-1})\f1{t(1-\ln(t))^2}dt.
\end{eqnarray}

Notice that
$$
\left\vert \int_{e^{-3}}^{e^{-3}-z-1}t^{-1}(1-\ln(t))^{-2}dt\right\vert \leq\f1 {16} e^3\theta
$$
and

\begin{eqnarray}
\nonumber&\ & \left\vert \int_{\theta}^{\theta-z-1}\f{t(1-\ln(t))^2}dt\right\vert \\
\nonumber&= & \left\vert (\ln(\theta-z-1)-\ln\theta)(1-\ln(\theta-z-1))^{-1}(1-\ln\theta)^{-1}\right\vert \\
\nonumber&\leq & \left\vert \frac{1}{2}\ln2+\frac{\pi}{4}i\right\vert (1-\frac{1}{2}\ln(2)-\ln\theta)^{-1}(1-\ln\theta)^{-1}.
\end{eqnarray}
By (\ref{i20}) and (\ref{i25}) we obtain
\begin{equation}\label{i2}
  I_2(z)\sim (1+z)^{-1}\ln^{-2}\theta\sim(1+z)^{-1}\ln^{-2}(1+z).
\end{equation}

Now by (\ref{dq}), (\ref{3parts}), (\ref{i12}), (\ref{i11}), (\ref{I3}) and (\ref{i2}), we obtain
\begin{equation}\label{ddq}
 \f d{dz}Q\tilde{a}(z)+(1+z)^{-1}(1-\ln\left\vert \theta\right\vert )^{-1}\sim\left\vert \theta\right\vert^{-1}\ln^{-2}\left\vert \theta\right\vert .
\end{equation}

It is easy to check that
\[
\textrm{Re}(iz(1+z)^{-1}(1-\ln\left\vert \theta\right\vert )^{-1})= \frac{\sin\theta(1-\ln \left\vert \theta\right\vert )^{-1}}{2-2\cos\theta}
\]
and
$$
\textrm{Im}(iz(1+z)^{-1}(1-\ln\left\vert \theta\right\vert )^{-1})= \frac12(1-\ln \left\vert \theta\right\vert )^{-1}.
$$
Now by
\[
 \frac{\sin\theta(1-\ln \left\vert \theta\right\vert )^{-1}}{2-2\cos\theta}+\theta^{-1}\ln^{-1}\left\vert \theta\right\vert  \sim  \theta ^{-1}\ln^{-2}\left\vert \theta\right\vert,
\]
and (\ref{ddq}) we obtain
\begin{eqnarray}\label{dreq}
  \nonumber \f d{d\theta}(\textrm{Re}Q\tilde{a}(z))+\theta^{-1}\ln^{-1}\left\vert \theta\right\vert &=&\textrm{Re}(iz\f d{dz}Q\tilde{a}(z))+\theta^{-1}\ln^{-1}\left\vert \theta\right\vert \\
 &\sim& \theta ^{-1}\ln^{-2}\left\vert \theta\right\vert .
\end{eqnarray}
\vskip2mm
\noindent By the similar argument, we can obtain
\begin{equation}\label{dimq}
 \f d{d\theta}(\textrm{Im}Q\tilde{a}(z)) \sim \left\vert \theta\right\vert ^{-1}\ln^{-2}\left\vert \theta\right\vert,
\end{equation}
although we will give another proof next.
\vskip2mm
\noindent
By (\ref{q}), $\textrm{Im}Q\tilde{a}(-e^{i\theta})$ is continuous at $\theta=0$  since
\begin{eqnarray}
\nonumber \textrm{Im}Q\tilde{a}(-e^{i\theta})&=&\int^0_{-1} \textrm{Im}\f1{(t+e^{i\theta})(1-\ln(1+ t))}dt\\
\nonumber& =& 
{-1}(\int_{-1}^{\theta-1}+\int_{\theta-1}^{e^{-1}-1}+\int_{e^{-1}-1}^{0})\frac{\sin\theta(1-\ln(1+t))^{-1}}{1+2t\cos\theta+t^{2}}dt\\
\nonumber&\leq&\int_{-1}^{\theta-1}\frac{\sin\theta(1-\ln(1+t))^{-1}}{(t-\cos\theta)^2+1-\cos^2\theta}dt\\
\nonumber&\ &+(\int_{\theta-1}^{e^{-1}-1}
    +\int_{e^{-1}-1}^{0})\frac{\sin\theta(1-\ln(1+t))^{-1}}{(1+t)^2+2t(\cos\theta-1)}dt\\
\nonumber& \leq  & \frac{\theta(1-\ln\theta)^{-1}}{\sin\theta}+
+\sin\theta(\int_{\theta-1}^{e^{-1}-1}
    +\int_{e^{-1}-1}^{0})\frac{(1-\ln(1+t))^{-1}}{(1+t)^{2}}dt\\
\nonumber&=&\frac{\theta(1-\ln\theta)^{-1}}{\sin\theta}+\sin\theta\int_{\theta-1}^{e^{-1}-1}((\f1{(1+t)(1-\ln(1+t)})'-(\f1{1-\ln(1+t)})'dt\\
\nonumber&\ &+\sin\theta\int_{e^{-1}-1}^{0}\frac{(1-\ln(1+t))^{-1}}{(1+t)^{2}}dt\\
& = &\frac{\theta(1-\ln\theta)^{-1}}{\sin\theta}-\frac{\sin\theta(1-\ln\theta)^{-1}}{\theta}+\sin\theta(1-\ln \theta)^{-1}+M\sin\theta,
\end{eqnarray}
where $M$ is a constant.

Notice that
$$
\textrm{Im}(t-z)^{-1}=-\f{\partial}{\partial t} \arctan(\frac{t+\cos\theta}{\sin\theta})
$$
and
$$
\f{\partial^{2}}{\partial t\partial \theta}\arctan(\frac{t+\cos\theta}{\sin\theta})>0
$$
and
$$
\frac{(1-t)(1-\ln(1+t))^{-2}}{1+2t\cos\theta+t^{2}}\geq 0.
$$
 Again by $1+2t\cos\theta+t^2\geq(1+t)^2$, We have
\begin{eqnarray}
 \nonumber\f d{d\theta}\textrm{Im}Q\tilde{a}(z)
\nonumber&= & \int_{-1}^{0}-(\f{\partial^{2}}{\partial t\partial \theta}\arctan(\frac{t+\cos\theta}{\sin\theta}))(1-\ln(1+t))^{-1}dt\\
\nonumber&= & -\f{\partial}{\partial \theta} \arctan\frac{t+\cos\theta}{\sin\theta}|_{t=0}-\int_{-1}^{0}\frac{(1+t\cos\theta)(1-\ln(1+t))^{-2}}{(1+t)(1+2t\cos\theta+t^{2})}dt\\
\nonumber&= & 1-\f12\int_{-1}^{0}(\frac{1}{1+t}+\frac{1-t}{1+2t\cos\theta+t^{2}})(1-\ln(1+t))^{-2}dt\\
\nonumber&= & \f12-\f12(\int_{-1}^{\theta-1}+\int_{\theta-1}^{e^{-3}-1}+\int_{e^{-3}-1}^{0})\frac{(1-t)(1-\ln(1+t))^{-2}}{1+2t\cos\theta+t^{2}})dt\\
\nonumber&\geq &\f12-\f12\int_{-1}^{\theta-1}\f{2(1-\ln\theta)^2}{1-\cos^2\theta}dt+\f12\int_{\theta-1}^{e^{-3}-1}\frac{2(1-\ln(1+t))^{-2}}{(1+t)^{2}}dt\\
\nonumber&\ &-\f12\int_{e^{-3}-1}^{0}\frac{(1-t)(1-\ln(1+t))^{-2}}{1+2t\cos\theta+t^{2}}dt\\
\nonumber&= & \f12-\theta\sin^{-2}\theta(1-\ln\theta)^{-2}-(1+t)^{-1}(1-\ln(1+t))^{-2}|_{t=e^{-3}-1}^{\theta-1}\\
         &\  &-\f12\int_{e^{-3}-1}^{0}\frac{(1-t)(1-\ln(1+t))^{-2}}{1+2t\cos\theta+t^{2}}dt,
\end{eqnarray}
and then
\begin{equation}\label{diq}
  0\geq\f d{d\theta}\textrm{Im}Q\tilde{a}(z) \geq -M\left\vert x\right\vert ^{-1}\ln^{-2}\left\vert x\right\vert.
\end{equation}

In short, $0\geq dImQ\tilde{a}(z)/dx \geq -M\left\vert x\right\vert ^{-1}\ln^{-2}\left\vert x\right\vert$, the first inequality holds because $\frac{(1-t)(1-\ln(1+t))^{-2}}{1+2t\cos(x)+t^{2}}\geq 0$ uniformly converges to $\frac{(1-t)(1-\ln(1+t))^{-2}}{1+2t+t^{2}}$ on $t\in [t_{0},0]$ for every $t_{0}\in (-1,0)$, which integral is diverse thus $\int_{-1}^{0}\frac{(1-t)(1-\ln(1+t))^{-2}}{1+2t\cos(x)+t^{2}}dt\geq 1$ for $x$ near 0 enough.

Now we obtain the following theorem.
\bT $Qa$  is $BMO_{\log}$ on $\mathbb{T}$.\eT
\vskip2mm
\noindent Now we prove the main theorems.
\bT
$N(T_{a})=0$ and $R(T_{a})$ is dense in $H^{1}(\mathbb{T})$
\eT
\bpp
Assume that there exists $f(t)\in H^1(\mathbb{T})$ such that $\mathrm{T}_a f(t)=P(a(t)f(t))=0$. Then $h(t)=a(t)f(t)$ must be in the complement of $H^(\mathbb{T})$, that is, $h(t)$ is the boundary value of a holomorphic function in $\mathbb{C}\setminus\mathbb{D}$ satisfying $h(\infty)=0$.

 Denote by $g(t)=a^{+}f=a^{-}h$. Then $g(t)$ can be extended be be a holomorphic  for all $z\in\mathbb{T}\backslash\{-1\}$. We are not sure whether $g(t)$ belongs to $L^{1}(\mathbb{T})$ due to the singularity at $z=-1$. But  $(1+t)g(t)$ must be in $ L^{1}(\mathbb{T})$. Then  $(1+z)g(z)$ is a bounded entire function. Especially, $(1+z)a^+(z)f(z)$ is a constant function on $\mathbb{D}$, which gives that $f\equiv0$ since $(1+z)a^+(z)$ is not constant.

  Now assume  $g\in H^{\infty}(\mathbb{T})$. Let $\eta(t)=(P(a^{-}g))/a^{+}$. It is obvious that  $\eta(t)\in H^1(\mathbb{T})$.  We obtain
  \begin{eqnarray}
   \nonumber T_a\eta &=& P(P(a^-g)\f1{a^-}) \\
   \nonumber  &=&  P(P(a^-g)\f1{a^-})+P(Q(a^-g)\f1{a^-})\\
   \nonumber  &=& P((P+Q)(a^-g)\f1{a^-})=P(a^-g\f1{a^-})=g.
  \end{eqnarray}
  So $\eta(t)$ is the preimage of $g(t)$.  Then the proof follows since $H^{\infty}(\mathbb{T})$ is dense in $H^1(\mathbb{T})$.
\epp

\bT
$R(T_a) \neq H^{1}(\mathbb{T})$.
\eT
\bpp
Let $g(z)=(1+z)^{-1}(a^{+}(z))^{-1}(\ln(a^{+}(z))+2)^{-2}$ for $z\in \overline{\mathbb{D}}=\mathbb{D}\cup \mathbb{T}$. It is easy to check that $g(t)\in H^{1}(\mathbb{T})$. We show that  $g(t)$  has no preimage of $T_a$ in $H^{1}(\mathbb{T})$.

Otherwise, assume that $f\in H^{1}(\mathbb{T})$ satisfies $T_af=g$.

Suppose that
   \begin{equation}\label{plk}
   af= \frac{a^{+}}{a^{-}}f=T_af+h=g+h.
\end{equation}
By Plemelj decomposition theorem\cite{mp}, we have $h\in \overline{H_{0}^{1}(\mathbb{T})}$.  we  multiply $(1+t)a^{-}$ on both sides of (\ref{plk}). And get $Q((1+t)a^{-}(g+h))=Q((1+t)a^+ f)=0$. Be similar to (\ref{q}), we have
\begin{equation}\label{}
 ( Q((1+t)a^{-}g))(z)=-\int_{-1}^{0}(t-z)^{-1}(1+t)g(t)dt,
\end{equation}
where $z\in \mathbb{C}\setminus [-1,0]$. Denote by  $G(z)=-(Q((1+t)a^{-}g))(z)$.
\vskip 2mm
We assume that  $z\in \mathbb{T}$ near $-1$ from now on.
\vskip 2mm
According to Lebesgue controlled convergent theorem, we have
\begin{equation}\label{}
  \lim_{z\rightarrow{-1}}G(z)=G(-1)=\int_{-1}^{0}g(t)dt.
\end{equation}
Denote by $h(z)=\sum^{\infty}_{n=1} h_{-n}z^{-n}$. Notice that  $(1+z)a^{-}(z)h(z)=P((1+z)a^{-}(z)h(z))+Q((1+z)a^{-}(z)h(z))$ and  $G(z)=-Q((1+z)a^{-}(z)h(z))$.
we have $G(z)+(1+z)a^{-}h=P((1+z)a^{-}h)=h_{-1}$. Then we have
\begin{equation}\label{hhh}
h=\frac{h_{-1}-G(z)}{(1+z)a^{-}(z)}\notin L^{1}(\mathbb{T}).
\end{equation}
If $h_{-1}\neq G(-1)$,it is obviously $h\notin L^{1}(\mathbb{T})$,  which contradicts to the assumption.

In the following, we prove that $h$ is still not in $L^{1}(\mathbb{T})$ when $h_{-1}=G(-1)$.

To estimate $ G(z)-h_{-1}=G(z)-G(-1)$ when $z=-\exp(i\theta)$ tends to $-1$, we decompose it into three parts as follows,
\begin{eqnarray}\label{G1}
 \nonumber  G(z)-G(-1)&=& \left(\int_{-1}^{\theta-1}+\int^{e^{-3}-1}_{\theta-1}+\int_{e^{-3}-1}^0\right)((t-z)^{-1}(1+t)g(t)-g(t))dt \\
   &=& I_1+I_2+I_3.
\end{eqnarray}

 Recall that  $\left\vert t-z\right\vert \geq \sin\theta$. We obtain
\begin{equation}
 \nonumber \left\vert\int_{-1}^{\theta-1}(t-z)^{-1}(1+t)g(t)dt\right\vert \leq \theta^2\sin^{-1}\theta g(\theta-1)\leq C|\ln\theta|^{-1}|\ln(|\ln\theta||^{-2}
\end{equation}
and
\begin{equation}
 \nonumber \int_{-1}^{\theta-1}g(t)dt=(\ln(1-\ln\theta)+2)^{-1},
\end{equation}
which give
\begin{equation}\label{i111}
  |I_1+\int_{-1}^{\theta-1}g(t)dt=(\ln(1-\ln\theta)+2)^{-1}|\leq C|ln\theta|^{-1}|\ln(|\ln\theta||^{-2}.
\end{equation}
It is easy to check
\[
\int_{e^{-3}-1}^{0}\frac{g(t)dt}{t-z}=(1+z)^{-1}\int_{e^{-3}-1}^{0}(\frac{1+t}{t-z}-1)g(t)dt
\]
converges to
\[
A=\int_{e^{-3}-1}^{0}(1+t)^{-1}g(t)dt,
\]
i.e.,
\begin{equation}\label{i333}
  I_3=\int_{\exp(-3)-1}^{0}((t-z)^{-1}(1+t)g(t)-g(t))dt\sim A(1+z)\sim\theta.
\end{equation}

Now we deal with the integral on the interval $[\theta-1,e^{-3}-1]$.
For convenience, denote by $W_t(\zeta)=(1-\ln(t-\zeta))^{-1}(\ln(1-\ln(t-\zeta))+2)^{-2}$.  Then $g(t-\zeta-1)=(t-\zeta)^{-1}W_t(\zeta)$, $g(t)=(1+t)^{-1}W_t(-1)$. And we have
\begin{eqnarray}
\nonumber W_t(z)-W_t(-1)&= & -\int_{0}^{\theta}\zeta(t-\zeta)^{-1}(1-\ln(t-\zeta))^{-2}((\ln(1-\ln(t-\zeta))+2)^{-2}\\
&\ &+2(\ln(1-\ln(t-\zeta))+2)^{-3})ds,
\end{eqnarray}
where $\zeta=-e^{is}$. Notice still $\left\vert t-\zeta\right\vert \leq \exp(-2)$ when $-1\leq t\leq \exp(-3)-1$. Again recall that  $|t-z|>|1+t|$, $|t-\zeta|>|1+t|$ and $s^{-1}(1-\ln s)^{-2}$ is decreasing when $s\in(0,\exp(-1)$. We have
\begin{eqnarray}\label{c221}
\nonumber
&\ & \left\vert \int_{\theta-1}^{\exp(-3)-1}(g(t-1-z)-(t-z)^{-1}(1+t)g(t))dt\right\vert \\
\nonumber
&= & \left\vert \int_{\theta-1}^{\exp(-3)-1}(t-z)^{-1}(W_t(z)-W_t(-1))dt\right\vert \\\nonumber&\leq & (3^{-2}+2\cdot 3^{-3})\int_{\theta-1}^{\exp(-3)-1}(1+t)^{-1}\int_0^\theta|t-\zeta|^{-1}(1-\ln|t-\zeta|)^{-2}dt\\
\nonumber&\leq &  \f 5 {27}\int_{\theta-1}^{\exp(-3)-1}(1+t)^{-1}\theta(1+t)^{-1}(1-\ln(1+t))^{-2}dt\\
\nonumber&= &  \f 5 {27}\theta\int_{\theta-1}^{\exp(-3)-1}((-(1+t)^{-1}(1-\ln(1+t))^{-2})'-3(1+t)^{-2}(1-\ln(1+t))^{-3})dt\\
\nonumber&\leq &  \f 5 {27}\theta\int_{\theta-1}^{\exp(-3)-1}((-(1+t)^{-1}(1-\ln(1+t))^{-2})'dt\\
&= & \f 5 {27}(1-\ln\theta)^{-3}-\f5{432}e^3\theta.
\end{eqnarray}
Applying Cauchy integral theorem on analytic function $g(t-1)$ along the parallelogram with vertices $x$, $e^{-3} $,  $e^{-3}-1 -z$ and $x-z-1$ taken in order, we obtain
\begin{equation}\label{}
\nonumber \left(\int_{\theta}^{e^{-3}}+\int^{\theta-1-z}_{\exp(-3)-1-z}\right)g(t-1)dt
= -\left (\int^{\theta}_{\theta-1-z}+\int_{\exp(-3)}^{\exp(-3)-1-z}\right)g(t-1)dt
\end{equation}
or
\begin{equation}
\nonumber \int_{\theta-1}^{e^{-3}-1}(g(t-1-z)-g(t))dt
=  \left(\int_{e^{-3}}^{e^{-3}-1-z}-\int_{\theta}^{\theta-z-1}\right)g(t-1)dt.
\end{equation}
But
$$
\left\vert \int_{e^{-3}}^{e^{-3}-1-z}g(t-1)dt\right\vert \leq \f {e^{3}}{36}\theta
$$
 and
\begin{eqnarray}
\nonumber \left\vert \int_{x}^{\theta-1-z}g(t-1)dt\right\vert
&= & \left\vert\frac{\ln(1-\ln\theta)-\ln(1-\ln(\theta-1-z))}{(\ln(1-\ln(\theta-1-z))+2)(\ln(1-\ln\theta)+2)}\right\vert \\
\nonumber&\leq & \left\vert \frac{1}{2}\ln2+\frac{\pi}{4}i\right\vert (1-\frac{1}{2}\ln2-\ln\theta)^{-1}\bullet\\
\nonumber&\ &(\ln(1-\frac{1}{2}\ln2-\ln\theta)+2)^{-1}(\ln(1-\ln\theta)+2)^{-1}\\
\nonumber&\leq &C|\ln\theta|^{-2}|\ln|\ln\theta||^{-1}.
\end{eqnarray}
So we obtain
\begin{equation}\label{c222}
\left\vert\int_{\theta-1}^{e^{-3}-1}(g(t-1-z)-g(t))dt\right\vert
\leq C|\ln\theta|^{-2}|\ln|\ln\theta||^{-1}.
\end{equation}

By (\ref{c221}) and (\ref{c222}), we obtain
\begin{eqnarray}\label{i222}
 \nonumber|I_2|  &=& \left\vert \int_{\theta-1}^{\exp(-3)-1}((t-z)^{-1}(1+t)g(t)-g(t))dt\right\vert \\
 \nonumber&\leq& \left\vert \int_{\theta-1}^{\exp(-3)-1}(g(t-1-z)-(t-z)^{-1}(1+t)g(t))dt\right\vert\\
 &\ &+\left\vert\int_{\theta-1}^{e^{-3}-1}(g(t-1-z)-g(t))dt\right\vert \leq M|\ln\theta|^{-2}.
\end{eqnarray}

Now by (\ref{i111}), (\ref{i333} and (\ref{i222}), we obtain
\begin{equation}\label{GGG}
  \left\vert G(z)-G(-1)+(\ln(1-\ln(x))+2)^{-1}\right\vert \leq C|\ln\theta|^{-1}|\ln\ln\theta|^{-2}.
\end{equation}

Recall that $h=\frac{G(-1)-G(z)}{(1+z)a^{-}(z)}\notin L^{1}(\mathbb{T})$. By (\ref{GGG}) we obtain
\[
h(z)-\frac{(\ln(1-\ln\theta))+2)^{-1}}{(1+z)a^{-}(z)}\in L^{1}(U).
\]
However,
\[
\frac{(\ln(1-\ln\theta)+2)^{-1}}{(1+z)a^{-}(z)}\notin L^{1}(U).
\]
So we get that $h\notin L^{1}(\mathbb{T})$, which still contradicts to the assumption.

\epp

Now we state the final result.
\bT
$T_a$ is not Fredholm.
\eT

\end{document}